\documentclass[12pt,a4paper]{article}
\usepackage{amscd,amsmath,amssymb,amsfonts,times}
\usepackage{mathrsfs}
\usepackage[dvips]{color} 
\usepackage[all]{xy}
\numberwithin{equation}{section}
    \newtheorem{thm}{Theorem}[section]
    \newtheorem{lem}[thm]{Lemma}

    \newtheorem{rem}[thm]{Remark}
\newcommand{\qed}
{\mbox{}\nolinebreak$\square$\medbreak\par}
\newenvironment{pf}{\par\smallskip\noindent\emph{Proof.}}{\hfill\qed\par\smallskip}
\newenvironment{pf*}[1]{\par\smallskip\noindent\emph{#1.}}{\hfill\qed\par\smallskip}
\begin{document}
\title{A simple construction of indecomposable higher Chow cycles 
in elliptic surfaces}
\author{M. Asakura}
\date\empty
\maketitle

\def\can{\mathrm{can}}
\def\Ext{\mathrm{Ext}}
\def\ch{{\mathrm{ch}}}
\def\Coker{\mathrm{Coker}}
\def\crys{\mathrm{crys}}
\def\dlog{d{\mathrm{log}}}
\def\dR{{\mathrm{d\hspace{-0.2pt}R}}}            
\def\et{{\mathrm{\acute{e}t}}}  
\def\id{{\mathrm{id}}}              
\def\Image{{\mathrm{Im}}}        
\def\Im{{\mathrm{Im}}}        
\def\Hom{{\mathrm{Hom}}}  
\def\ker{{\mathrm{Ker}}}          
\def\mf{{\text{mapping fiber of}}}
\def\Pic{{\mathrm{Pic}}}
\def\CH{{\mathrm{CH}}}
\def\NS{{\mathrm{NS}}}
\def\NF{{\mathrm{NF}}}
\def\End{{\mathrm{End}}}
\def\pr{{\mathrm{pr}}}
\def\Proj{{\mathrm{Proj}}}
\def\ord{{\mathrm{ord}}}
\def\qis{{\mathrm{qis}}}
\def\reg{{\mathrm{reg}}}          %
\def\res{{\mathrm{res}}}          %
\def\Res{\mathrm{Res}}
\def\Spec{{\mathrm{Spec}}}     
\def\MHS{{\mathrm{MHS}}}
\def\cont{{\mathrm{cont}}}
\def\rank{{\mathrm{rank}}}
\def\sp{{\mathrm{sp}}}

\def\bA{{\mathbb A}}
\def\bC{{\mathbb C}}
\def\C{{\mathbb C}}
\def\G{{\mathbb G}}
\def\bE{{\mathbb E}}
\def\bF{{\mathbb F}}
\def\F{{\mathbb F}}
\def\bG{{\mathbb G}}
\def\bH{{\mathbb H}}
\def\bJ{{\mathbb J}}
\def\bL{{\mathbb L}}
\def\cL{{\mathscr L}}
\def\bN{{\mathbb N}}
\def\bP{{\mathbb P}}
\def\P{{\mathbb P}}
\def\bQ{{\mathbb Q}}
\def\Q{{\mathbb Q}}
\def\bR{{\mathbb R}}
\def\R{{\mathbb R}}
\def\bZ{{\mathbb Z}}
\def\Z{{\mathbb Z}}
\def\cH{{\mathscr H}}
\def\cD{{\mathscr D}}
\def\cE{{\mathscr E}}
\def\cO{{\mathscr O}}
\def\O{{\mathscr O}}
\def\cR{{\mathscr R}}
\def\cS{{\mathscr S}}
\def\cX{{\mathscr X}}
\def\cZ{{\mathscr Z}}
\def\cC{{\mathscr C}}
%
\def\ep{\epsilon}
\def\vG{\varGamma}
\def\vg{\varGamma}
\def\cyc{{\mathrm cycl}}
%
%
%
\def\lra{\longrightarrow}
\def\lla{\longleftarrow}
\def\Lra{\Longrightarrow}
\def\hra{\hookrightarrow}
\def\lmt{\longmapsto}
\def\ot{\otimes}
\def\op{\oplus}
\def\wt#1{\widetilde{#1}}
\def\wh#1{\widehat{#1}}
\def\spt{\sptilde}
\def\ol#1{\overline{#1}}
\def\ul#1{\underline{#1}}
\def\us#1#2{\underset{#1}{#2}}
\def\os#1#2{\overset{#1}{#2}}
\def\lim#1{\us{#1}{\varinjlim}}
\def\plim#1{\us{#1}{\varprojlim}}

\begin{center}
{\it Dedicated to the 60th birthday of Professor James D. Lewis}
\end{center}

\section{Introduction}
Let $\CH^j(X,i)$ be Bloch's
higher Chow groups of a projective smooth variety $X$ over $\C$.
A higher Chow cycle $z\in \CH^j(X,i)$ is called {\it indecomposable}
if it does not belong to the image of the map of the product
\[
\CH^1(X,1)\ot \CH^{j-1}(X,i-1)\lra \CH^j(X,i).
\]
Of particular interest to us is $\CH^2(X,1)$.
For $A=\Q$ or $\R$,
we say $z\in \CH^2(X,1)$ {\it $A$-regulator indecomposable}
if the regulator class $\reg(z)\in H^3_\cD(X,A(2))$
in the Deligne-Beilinson cohomology group with coefficients in $A$
does not belong to the image of $H^1_\cD(X,\Z(1))\ot H^2_\cD(X,A(1))\cong 
\C^\times\ot \CH^1(X)\ot A$.
In other words, $z$ is $A$-regulator indecomposable if and only if
\[
\reg(z)\ne0 \in \mathrm{Ext}_{\mathrm{MHS}}^1(A,H^2_{\mathrm{ind}}(X,A(2))),
\quad
H^2_{\mathrm{ind}}(X,A):=H^2(X,A)/\NS(X)\ot A.
\]
Obviously
$
\R\mbox{-reg. indecomp.}\Longrightarrow
\Q\mbox{-reg. indecomp.}\Longrightarrow
\mbox{indecomposable}.
$

\medskip

Quite a lot of examples of $\Q$ or $\R$-regulator indecomposable cycles are obtained by many
people (\cite{AS-real}, \cite{lewis}, \cite{lewisJAG}, \cite{kerr-indk1}, \cite{muller-stach}
and more).

In this note we construct $\R$-regulator indecomposable cycles
for $X$ an elliptic surface which satisfies certain conditions.
The main theorem is the following.

\begin{thm}\label{main}
Let $S$ be a smooth irreducible curve over $\C$. Let 
\[
\xymatrix{
\cX\ar[rd]_g\ar[rr]_f&&\cC\ar[ld]^h\ar@/_/[ll]_s\\
&S&
}
\]
be an elliptic fibration over $S$ with a section $s$. 
This means that $g$ and $h$ are
projective smooth morphisms of relative dimension $2$ and $1$ respectively, and
the general fiber of $f$ is an elliptic curve. 
For a point $t\in S$ we denote $X_t=g^{-1}(t)$ or $C_t
=h^{-1}(t)$ the fibers
over $t$. 
Assume that the following conditions hold.
\begin{enumerate}
\item[(1)]
Let $\eta$ be the generic point of $S$. 
Then there is a split multiplicative fiber $D_\eta=f^{-1}(P)\subset X_\eta$ of Kodaira
type $I_n$, $n\geq 1$.
\item[(2)]
Let $\cD\subset \cX$ be the closure of $D_\eta$. Then
there is a closed point $0\in S(\C)$ such that the specialization 
$D_0:=\cD\times_\cX X_0$ is multiplicative
of type $I_m$ with $m>n$.
\end{enumerate}
Then the composition
\[
\CH^1(D_t,1)\lra \CH^2(X_t,1)\os{\reg}{\lra} \Ext^1_\MHS(\R,H^2(X_t,\R(2))/\NF(X_t)\ot\R)
\]
is non-zero for a general $t\in S(\C)$.
Here $\NF(X_t)\subset \NS(X_t)$ denotes the subgroup generated by components
of singular fibers and the section $s(C_t)$.
In particular, if $\NF(X_t)\ot\Q=\NS(X_t)\ot\Q$, then there is a $\R$-regulator
indecomposable higher Chow cycle.
\end{thm}
The key assumption is a ``degeneration of $D_t$", which often appears in
a family of elliptic surfaces.

\bigskip

\begin{center}
\unitlength 0.1in
\begin{picture}( 41.3100, 21.7000)( 13.2900,-25.7000)
%
{\color[named]{Black}{%
\special{pn 8}%
\special{pa 4780 2470}%
\special{pa 3810 1224}%
\special{fp}%
}}%
%
{\color[named]{Black}{%
\special{pn 8}%
\special{ar 2120 1340 650 900  1.5378593  4.7123890}%
}}%
%
{\color[named]{Black}{%
\special{pn 8}%
\special{ar 1350 1320 650 900  4.7123890  1.6037333}%
}}%
\put(16.0000,-27.0000){\makebox(0,0)[lb]{$D_t$}}%
\put(28.0000,-19.3000){\makebox(0,0)[lb]{$t\lra 0$}}%
%
{\color[named]{Black}{%
\special{pn 8}%
\special{pa 4780 400}%
\special{pa 3810 1648}%
\special{fp}%
}}%
%
{\color[named]{Black}{%
\special{pn 8}%
\special{pa 4490 2470}%
\special{pa 5460 1224}%
\special{fp}%
}}%
%
{\color[named]{Black}{%
\special{pn 8}%
\special{pa 4490 400}%
\special{pa 5460 1648}%
\special{fp}%
}}%
%
{\color[named]{Black}{%
\special{pn 8}%
\special{pa 2690 1450}%
\special{pa 3430 1450}%
\special{fp}%
\special{sh 1}%
\special{pa 3430 1450}%
\special{pa 3364 1430}%
\special{pa 3378 1450}%
\special{pa 3364 1470}%
\special{pa 3430 1450}%
\special{fp}%
}}%
\put(45.6000,-27.0000){\makebox(0,0)[lb]{$D_0$}}%
\end{picture}%
\end{center}
\begin{center}
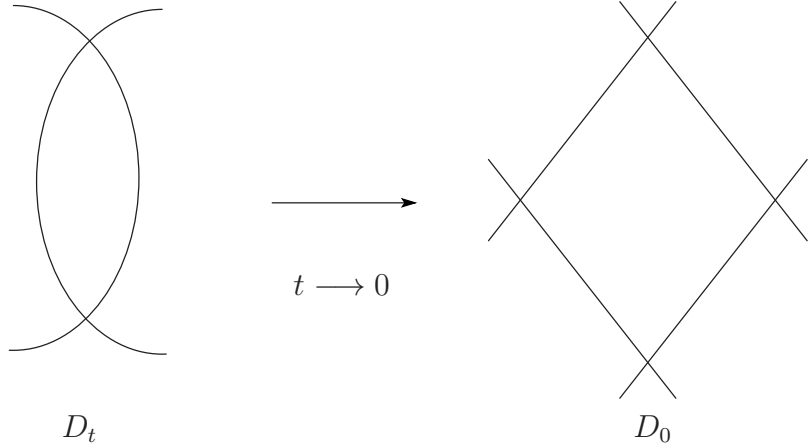
Figure : Degeneration of $I_2$
to $I_4$
\end{center}

\bigskip

In  \S \ref{GL-sect}, we will apply Theorem \ref{main} to construct a
$\R$-reg. indecomp. cycle in a self-product of elliptic curves.
However, to do it in more general situation,
the computation of the Picard number might be an obstacle. 
Indeed it is easy to compute the rank of $\NF(X_t)$, 
whereas there is no general method to do it for $\NS(X_t)$, and usually
it is done by case-by-case analysis.

\medskip

\noindent{\bf Acknowledgment.}
The author is grateful to Professor James D. Lewis for reading the first draft carefully
and giving many comments.

\section{Proof of Theorem \ref{main}}
We keep the notation and assumption in Theorem \ref{main}. 
\subsection{Step 1 : Construction of a higher Chow cycle}
\begin{lem}\label{km3}
Let $f:X\to C$ be an elliptic fibration over a field $K$ of characteristic$\ne 2,3$, 
with a (fixed) section $s$. 
Let $D=f^{-1}(P)$ be a split multiplicative fiber, and let $i:D\hra X$.
Then there is an exact sequence
\begin{equation}\label{km7}
\CH^1(\wt{D},1)_\Q\lra \CH^1(D,1)_\Q\os{v}{\lra} \Q\lra 0
\end{equation}
where $\wt{D}\to D$ is the normalization.
There exists a higher Chow cycle $Z\in \CH^2(D,1)$ such that $v(Z)\ne0$
and $i_*(Z)\in \CH^2(X,1)$ is vertical to generators of $\NF(X)$. 
Here ``vertical to $E$" means that it lies in the kernel of
the composition
$\CH^2(X,1)\to \CH^2(\wt{E}\times_K \ol{K},1)\to \ol{K}^\times\ot\Q$.
\end{lem}
\begin{pf}
We omit to show the exact sequence \eqref{km7} (easy exercise). 
We show the existence of $Z$.
In this proof, we use $K$-groups rather than higher Chow groups.
Let $D=\sum_i D_i$ be the irreducible decomposition. 
It is enough to construct $Z\in K'_1(D)^{(1)}$ such that $v(Z)\ne0$ and 
$i_*(Z)$ is vertical to each $D_i$ and $s(C)$ because it is obviously vertical to
the other fibral divisors.
Since $D$ is split multiplicative,
each $D_i$ is geometrically irreducible and the singularities of $D$ are $K$-rational.
Therefore we may assume $K=\ol{K}$ by the standard norm argument.

It is enough to show that
the image of the composition
\[
K'_1(D)^{(1)}\os{i_*}{\lra} K_1(X)^{(2)}\os{i^*}{\lra} 
K_1(\wt{D})^{(2)}\cong (K^\times\ot\Q)^{\op n}
\]
coincides with that of
\[
K'_1(\wt{D})^{(1)}\os{i_*}{\lra} 
K_1(X)^{(2)}\os{i^*}{\lra} K_1(\wt{D})^{(2)}\cong (K^\times\ot\Q)^{\op n}.
\]
Indeed, the above implies that there is a cycle $Z_0\in K'_1(D)^{(1)}$ such that
$v(Z_0)\ne0$ and $i_*(Z_0)$ is vertical to each $D_i$.
Let $f^*:K_1(K)\to K'_1(D)^{(1)}$. Then $Z:=Z_0+f^*(\lambda)$ for suitable 
$\lambda\in K^\times$ can be vertical to the section $s(C)$ and $D_i$.
Moreover $v(Z)=v(Z_0)\ne0$.

To do the above we may replace $X$ with 
$\hat{X}=f^{-1}(\Spec K[[s]])$ the formal neighborhood around $D$.
Then it is enough to show that there is $Z_0\in K'_1(D)^{(1)}$ such that
$v(Z_0)\ne0$ and $i_*(Z_0)=0$ in $K_1(\hat{X})^{(2)}$.
To do this, we may further replace $\hat{X}$ with 
$\hat{X}_n$ the minimal desingularization of
$\hat{X}\times\Spec K[[s^{1/n}]]$ for some $n\geq 1$
due to a commutative diagram
\[
\xymatrix{
\Q\ar[d]_n&K'_1(D_n)^{(1)}\ar[r]^{i_*}\ar[d]^{\phi_*}\ar[l]_{v\quad}
&K_1(\hat{X}_n)^{(2)}\ar[d]^{\phi_*}\\
\Q&K'_1(D)^{(1)}\ar[r]^{i_*}\ar[l]_{v\quad}&K_1(\hat{X})^{(2)}
}
\]
where $\phi:\hat{X}_n\lra \hat{X}$.
Thus we can assume $\hat{X}$ is defined by a Weierstrass equation
\[
y^2=x^3+x^2+c(s), \quad c(s)\in sK[[s]].
\]
Then letting $\partial:K_2(\hat{X}\setminus D)^{(2)}\to K'_1(D)^{(1)}$
be the boundary map, we put
\[
Z_0:=\partial\left\{
\frac{y-x}{y+x},\frac{-c}{x^3}
\right\}.
\]
This satisfies $v(Z_0)\ne0$ and $i_*(Z_0)=0$ in $K_1(\hat{X})^{(2)}$.
\end{pf}

By Lemma \ref{km3}, there is a higher Chow cycle $Z_\eta\in \CH^1(D_\eta,1)$ such that
$v(Z_\eta)\ne0$ and it is vertical to $\NF(X_\eta)$.
We use the same symbol ``$Z_\eta$" for $i_*(Z_\eta)\in \CH^2(X_\eta,1)$ with $i:D_\eta
\hra X_\eta$ 
since it will be clear from the context which is meant.

Let $\cZ^*\subset\cX^*$ be the closure of $Z_\eta$ in $\cX^*:=g^{-1}(S^*)$ 
for some nonempty Zariski open $S^*\subset S$.
Then the goal is to show nonvanishing
\begin{equation}\label{km1}
\reg(Z_t)\ne 0\in \Ext^1_\MHS(\R,H^2(X_t,\R(2))), \quad Z_t:=\cZ^*|_{X_t}
\end{equation}
for a general $t\in S^*(\C)$.
Indeed, since $Z_t$ is vertical to $\NF(X_t)$, the above implies
the desired nonvanishing
\[
\reg(Z_t)\ne 0\in \Ext^1_\MHS(\R,H^2(X_t,\R(2))/\NF(X_t)\ot\R).
\]
\subsection{Step 2 : Boundary of $Z_\eta$}
\begin{lem}\label{km2}
Let $\partial:\CH^2(X_\eta,1)\to \CH^1(X_0)$ be the boundary map arising from
the localization exact sequence.
Then $\partial(Z_\eta)$ is non-torsion.
\end{lem}
\begin{pf}
We may assume $S=\Spec R$ where $R$ is a DVR, with a closed point $0$ and
generic point $\eta$.
Let $D'_\eta\subset D_\eta$ and $D'_0\subset D_0$
be the unique chain of rational curves which forms
Neron polygons. Then there are exactly $n$ reduced components 
(resp. $m$ reduced components) in $D'_\eta$ (resp. $D'_0$), as it is of Kodaira type
$I_n$ (resp. $I_m$).
Since $m>n$, there is a reduced irreducible component $E'_\eta\subset D'_\eta$ such that
its specialization $E'_0$ has at least two reduced components.
Let $\cE'\subset \cD$ be the closure of $E'_\eta$ and $j:\wt{\cE}\to \cE'$ the normalization.
Let $\wt{E}_0:=\wt{\cE}\times_S\{0\}=\sum_{j=1}^q r_jC_j$ be the special fiber.
\[
\xymatrix{
\wt{E}_\eta\ar[r]\ar[d]&\wt{\cE}\ar[d]_j&\wt{E}_0\ar[d]\ar[l]\\
E'_\eta\ar[r]\ar[d]& \cE'\ar[d]\ar[d]&E'_0\ar[d]\ar[l]\\
X_\eta\ar[r]&\cX&X_0\ar[l]
}
\]
The generic fiber $\wt{E}_\eta:=\wt{\cE}\times_S\{\eta\}$ 
is a smooth irreducible rational curve.
Let $T_1,T_2\subset \wt{\cE}$ be the inverse image of the intersection locus
of the Neron polygon $D'_\eta$.
Let $T_1$ hits a component of $\wt{E}_0$, say $C_1$.
Then $T_2$ hits another component (say $C_2$), since
the image $j(\wt{E}_0)\subset D'_0$ has at least two component.
One has $r_1=r_2=1$, namely $C_1$ and $C_2$ are reduced components,
since $T_1$ and $T_2$ are sections of $\wt{\cE}\to S$.
Moreover the intersection points $(T_1\cap C_1)$ and $(T_2\cap C_2)$ are nonsingular points
in $\wt{\cE}$.

Let $f$ be a rational function on $E_\eta$ such that the divisor
$\mathrm{div}_{E_\eta}(f)=T_1-T_2$ and hence
\begin{equation}\label{bd1p}
\mathrm{div}_{\wt{\cE}}(f)=T_1-T_2+\sum_{j=1}^q a_jC_j.
\end{equation}
Then
\[
\partial(Z_\eta)=\sum_{j=1}^q a_jC_j+N\cdot j(\wt{E}_0)+(\mbox{other components of }D_0)
\in\CH^1(X_0)
\]
for some $N\in \Z$ by the construction of $Z_\eta$. 
Therefore it is enough to show $a_1\ne a_2$.

We take a desingularization $\rho:\cE\to \wt{\cE}$ (if necessary).
Then 
\begin{equation}\label{bd1}
\mathrm{div}_{\cE}(f)=T_1-T_2+\sum_{j=1}^q a_jC'_j+\sum_{j=q+1}^\ell a_jC'_j.
\end{equation}
where $C'_j$ ($j\leq q$) are the strict transform of $C_j$ and
$C'_j$ ($j>q$) are the exceptional curves.
We may assume $T_i$ intersects with 
$C'_i$ in $\cE$ for $i=1,2$, since $T_i\cap C_i$ are nonsingular points of $\wt{\cE}$.
Renumbering $C'_j$, we may assume $T_1$ intersects with $C'_1$ and
$T_2$ intersects with $C'_\ell$. We then want to show $a_1\ne a_\ell$.
By the intersection theory, \eqref{bd1} yields
\[\begin{cases}
1+\sum_{j=1}^\ell a_j(C'_1,C'_j)=0\\
\sum_{j=1}^\ell a_j(C'_i,C'_j)=0&2\leq \forall i\leq \ell-1
\end{cases}
\]
namely
\begin{equation}\label{bd2}
\begin{pmatrix}
(C'_1,C'_1)&\cdots&(C'_1,C'_\ell)\\
\vdots&&\vdots\\
(C'_{\ell-1},C'_1)&\cdots&(C'_{\ell-1},C'_\ell)
\end{pmatrix}
\begin{pmatrix}
a_1\\
\vdots\\
a_\ell
\end{pmatrix}
=\begin{pmatrix}
-1\\
\vdots\\
0
\end{pmatrix}
\end{equation}
By replacing $(a_1,\cdots,a_\ell)$ with
$(a_1,\cdots,a_\ell)+c(r_1,\cdots,r_\ell)$ for some $c\in \Q$, we may assume
$a_\ell=0$.
Then
\begin{equation}\label{bd3}
A
\begin{pmatrix}
a_1\\
\vdots\\
a_{\ell-1}
\end{pmatrix}
=
\begin{pmatrix}
(C'_1,C'_1)&\cdots&(C'_1,C'_{\ell-1})\\
\vdots&&\vdots\\
(C'_{\ell-1},C'_1)&\cdots&(C'_{\ell-1},C'_{\ell-1})
\end{pmatrix}
\begin{pmatrix}
a_1\\
\vdots\\
a_{\ell-1}
\end{pmatrix}
=\begin{pmatrix}
-1\\
\vdots\\
0
\end{pmatrix}
\end{equation}
By Zariski's lemma (\cite{barth} III (8.2)), one has $\det A<0$ and $\det A_{11}<0$ where
$A_{11}$ is the cofactor matrix.
Therefore one has $a_1=-\det A_{11}/\det A\ne 0$, the desired assertion. 
This completes the proof.
\end{pf}
\subsection{Step 3 : Extension of admissible variations of MHS's}
Let
\[
\reg(\cZ^*)\in H^3_{\cD}(\cX^*,\Q(2))
\]
be the Deligne-Beilinson cohomology class.
Lemma \ref{km2} together with the commutative diagram
\begin{equation}\label{km10}
\xymatrix{
\CH^2(\cX^*,1)\ar[r]^\partial\ar[d]_\reg&\CH^1(X_0)\ar[d]^{cl}\\
H^3_{\cD}(\cX^*,\Q(2))\ar[r]^{\Res}&H^2_{\cD}(X_0,\Q(1))
\ar[r]^{\cong\quad}
& H^2(X_0,\Q(1))\cap H^{1,1}
}
\end{equation}
yields nonvanishing $\reg(\cZ^*)\ne 0$.
There is the isomorphism
\begin{equation}\label{km15}
H^3_{\cD}(\cX^*,\Q(2))\cong 
\Ext_{\mathrm{VMHS}(S^*)}^1(\Q,H_\Q),\quad H_\Q:=R^2\pi_*\Q(2)
\end{equation}
with the extension group of 
admissible variations of mixed Hodge structures on $S^*$.
Hence
$\reg(\cZ^*)$ defines a non-trivial extension
\begin{equation}\label{km11}
0\lra H_\Q\lra V_\Q\lra \Q\lra0
\end{equation}
of admissible VMHS's on $S^*$.
Tensoring with $\R$, one has an extension
\begin{equation}\label{km11real}
0\lra H_\R\lra V_\R\lra \R\lra0
\end{equation}
of real VMHS's. This is also a non-trivial extension.
Then the following lemma finishes the proof of Theorem \ref{lewis}.
\begin{lem}\label{limit-lem}
Let $\Delta\subset S$ be a small neighborhood of $0$, and put 
$\Delta^*=\Delta\setminus\{0\}$.
For $t\in\Delta^*$ such that $0<|t|\ll 1$, the extension 
\begin{equation}\label{km11fiber}
0\lra H_{\R,t}\lra V_{\R,t}\lra \R\lra0
\end{equation}
is non-trivial
where $H_{\R,t}$ etc. denotes the fiber at $t$.
Hence the nonvanishing \eqref{km1} follows.
\end{lem}
Before proving Lemma \ref{limit-lem}, we note that
``$\Q$-regulator indecomposability of $Z_t$" is immediate from the fact
that \eqref{km11} is non-trivial.
Let $H_\O:=\O_{S^*}\ot H_\Q$ and $F^\bullet$ be the Hodge bundles.
Put $J_\Q:=H_\O/(F^2+H_\Q)$
and 
\[
J^h_\Q:=\ker[J_\Q\os{\nabla}{\lra} \Omega^1_{S^*}\ot H_\O/F^1]
\]
the sheaf of horizontal sections where $\nabla$ denotes the Gauss-Manin connection.
As is well-known, there is the injective map
\begin{equation}\label{km18}
\iota:\Ext_{\mathrm{VMHS}(S^*)}^1(\Q,H_\Q)\hra \vg(S^*,J^h_\Q),
\end{equation}
and $\nu_{\cZ^*}:=\iota(\reg(\cZ^*))\in\vg(S^*,J^h_\Q)$ is called
the {\it normal function}
associated to $\cZ^*$.
Since the zero locus of the normal function  
is at most a countable set, we have
$\nu_{\cZ^*}(t)\ne0$ for a general $t$, and hence
that $Z_t$ is $\Q$-regulator indecomposable.
However to obtain the ``$\R$-regulator indecomposability" in the same way, 
we need to show that the injectivity
of \eqref{km18} remains true if we replace $J^h_\Q$ with $J^h_\R$, and I don't know how to 
prove it in general. We prove Lemma \ref{limit-lem} in a different way.

\medskip

\noindent{\it Proof of Lemma \ref{limit-lem}.}

Let $T$ be the local monodromy around $\Delta^*$. 
The action on $H_t=H^2(X_t,\Q(2))$ is trivial, whereas
that on $V_{\Q,t}$ is non-trivial.
Indeed, let $\{e_{1,t},\cdots,e_{m,t}\}$ be a basis of $H_{\Q,t}$ and 
$\{e_{0,t},e_{i,t}\}$ be a basis of $V_{\Q,t}$ such that $e_{0,t}\not\in H_{\Q,t}$.
Let $N=T-1$ be the log monodromy.
Then one has $(2\pi i)^{-1}Ne_{0,t}=\Res(\reg(\cZ^*))=cl[\partial\cZ^*]\ne 0$ 
under the natural isomorphism 
$H_{\Q,t}\ot\Q(-1)\cong H_{\Q,0}\ot\Q(-1)=H^2(X_0,\Q(1))$.

Let us fix a frame $z_t\in H_{\Q,t}$ which satisfies that $z|_{t=0}
=cl[\partial\cZ^*]\in H^2(X_0,\C)$. Then the {\it admissibility}
of $V$ yields that there are
holomorphic functions $u_i(t)$ on $\Delta^*$ with at most meromorphic 
singularities at $t=0$ such that
\[
f_0=e_{0,t}+\frac{\log(t)}{2\pi i} z_t+\sum_i u_i(t)e_{i,t}
\]
belongs to the Hodge bundle $\vg(\Delta^*,F^0V_\O)$.
Moreover, by adding some $\theta\in \vg(\Delta^*,F^2H_\O)$, one has
\[
f'_0=f_0+\theta=e_{0,t} +\frac{\log(t)}{2\pi i} z_t+\sum_i u'_i(t)e_{i,t}
\]
with $u'_i(t)$ holomorphic at $t=0$.
Therefore the extension data of \eqref{km11} is given as follows
\begin{equation}\label{km12}
\mbox{(extension class of }V|_{\Delta^*})=\frac{\log(t)}{2\pi i} z_t+\sum_i u'_i(t)e_{i,t}
\in H_\O/(F^2+H_\Q).
\end{equation}
Let $\mathscr{A}_{\Delta^*}$ be the sheaf of $C^\infty$-functions on $\Delta^*$,
and let $c:H_\R\ot_\R\mathscr{A}_{\Delta^*}\to H_\R\ot_\R\mathscr{A}_{\Delta^*}$
the complex conjugation given by $c(x\ot f)=x\ot\bar{f}$.
Let 
\begin{equation}\label{km13}
H_\O/(F^2+H_\R)\os{1-c}{\lra} H_\R\ot\mathscr{A}_{\Delta^*}/(1-c)F^2\cong
(H_\R\cap H^{1,1})\ot\mathscr{A}_{\Delta^*}.
\end{equation}
Then the extension class \eqref{km12} goes to
\begin{equation}\label{km14}
\frac{\log|t|}{\pi} z_t+\sum_i \mbox{(bounded funcion)}e_{i,t}\in
(H_\R\cap H^{1,1})\ot\mathscr{A}_{\Delta^*}
\end{equation}
via \eqref{km13}. This does not vanish
for $0<|t|\ll 1$, and hence the proof is done. \qed

\medskip

This completes the proof of Theorem \ref{main}.

\section{Example : Self-Product of elliptic curves}\label{GL-sect}
Recall the following theorem due to Lewis and Gordon
\footnote{Though there was an error in their proof, 
T\"urkmen recently corrected it \cite{turkmen}.}.
\begin{thm}[Lewis - Gordon, \cite{lewisJAG} Thm.1]\label{lewis}
Let $X=E\times E$ be a product of elliptic curves over $\C$.
If $E$ is sufficiently general, then there exists a $\R$-regulator indecomposable cycle
$z\in \CH^2(X,1)$
\footnote{
\cite{lewisJAG} Thm.1 deals with a general product $E_1\times E_2$, though
their techniques allow to handle the case of a self-product $E\times E$ as well.}.
\end{thm}

We here give an alternative proof as an application of Theorem \ref{main}.
Let
\[
E_a:y^2=F(x)=x(x-1)(x-a), \quad a\in\C\setminus\{0,1\}
\]
be an elliptic curve over $\C$.
Let $X_a:=\mathrm{Km}(E_a\times E_a)$ be the {\it Kummer surface} associated to the product
$E_a\times E_a$.
The birational model of $X_a$ is given by the double sextic 
\[
X_a:w^2=F(x_1)F(x_2)=x_1x_2(x_1-1)(x_2-1)(x_1-a)(x_2-a).
\]
Changing the variables $x:=x_1$, $t:=x_1x_2$ and $y:=wx_1$, we get 
\begin{equation}\label{e-km1}
X_a:y^2=t(x-1)(x-t)(x-a)(ax-t).
\end{equation}
This gives an elliptic fibration
\begin{equation}\label{e-km2}
f:X_a\lra \P^1,\quad (x,y,t)\longmapsto t
\end{equation}
where $t$ is the affine parameter of $\P^1$.
This is naturally extended to a family over $S=\bA_\C\setminus\{0,1\}$
\[
\xymatrix{
\cX\ar[rd]_\pi\ar[rr]^f&&\P^1\times S\ar[ld]\\
&S&}
\]
in which $X_a=\pi^{-1}(a)$.
We fix a section $s$ of $(x,y)=(1,0)$.

Let us look at the singular fibers in \eqref{e-km2}.
The singular fibers are located at $t=1,a,a^2,0,\infty$.
Two additive fibers appear at $t=0,\infty$ and both are of Kodaira type $I_2^*$.
If $a\ne -1$, the singular fibers at $t=1,a,a^2$ are multiplicative
of type $I_2$, $I_4$, $I_2$ respectively, and
if $a=-1$, two fibers at $t=\pm1$ are multiplicative
of type $I_4$. 
Let $\NF(X_a)\subset \NS(X_a)$ be the subgroup generated by irreducible components of
singular fibers and the section $s$. 
The rank of $\NF(X_a)$ is $19$ if $a\ne-1$, and $20$ if $a=-1$.
As is well-known, 
\begin{equation}\label{cm0}
\rank \NS(X_a)=
\begin{cases}
19&E_a\mbox{ has no CM}\\
20&E_a\mbox{ has a CM}.
\end{cases}
\end{equation}
In particular 
\begin{equation}\label{e-km3}
\NF(X_a)\ot\Q= \NS(X_a)\ot\Q \mbox{ if and only if $E_a$ has no CM.}
\end{equation}
If $a\ne-1$, there is the isomorphism
\begin{equation}\label{e-km4}
H^2(X_a,\Q)/\NF(X_a)_\Q\cong \mathrm{Sym}^2 H^1(E_a,\Q).
\end{equation}

We now apply Theorem \ref{main} for $\cD=f^{-1}(1)$ and $X_{-1}=\pi^{-1}(-1)$.
As we see in the above,
$D_\eta$ is a split multiplicative fiber of type $I_2$ and
$D_{-1}:=\cD\times_\cX X_{-1}$ is of type $I_4$, so the conditions (1) and (2) are
satisfied.
We obtain a higher Chow cycle $Z_a\in \CH^2(X_a,1)$ arising from $f^{-1}(1)$,
and this is $\R$-regulator indecomposable for a general $a$
by Thm. \ref{main} and \eqref{e-km3}. By \eqref{e-km4}, this gives
a $\R$-regulator indecomposable cycle in a self product $E_a\times E_a$. 

\begin{rem}
The real regulator $\reg(Z_a)$ (as a function of $a$) is studied 
in detail in \cite{lewis} \S 6.4.
\end{rem}


\bigskip

\noindent
Department of Mathematics, Hokkaido University,
Sapporo 060-0810,
JAPAN

\medskip

\noindent
asakura@math.sci.hokudai.ac.jp

\end{document}